\documentclass[a4paper,reqno,11pt,oneside]{amsart}\usepackage[left=3cm,right=3cm,top=3.5cm,bottom=3.5cm]{geometry}
\usepackage[colorlinks=true,urlcolor=blue,
citecolor=red,linkcolor=blue,linktocpage,pdfpagelabels,
bookmarksnumbered,bookmarksopen]{hyperref}
\usepackage[english]{babel}
\usepackage{graphicx}
\usepackage[toc,page,header]{appendix}
\usepackage{mathrsfs}
\usepackage{amssymb,amsmath, amsfonts, amsthm, mathtools}
\usepackage{xcolor}

\makeatletter
\@namedef{subjclassname@2020}{%
  \textup{2020} Mathematics Subject Classification}
\makeatother

\usepackage{latexsym}

\allowdisplaybreaks[1]

\numberwithin{equation}{section}

\DeclareMathOperator{\divergence}{div}

\renewcommand{\(}{\left(}
\renewcommand{\)}{\right)}
\renewcommand{\[}{\left[}
\renewcommand{\]}{\right]}

\newtheorem{theorem}{Theorem}[section]
\newtheorem{proposition}[theorem]{Proposition}

\newtheorem{lemma}[theorem]{Lemma}
\theoremstyle{definition}
\newtheorem{remark}[theorem]{Remark}
\theoremstyle{definition}
\newtheorem{definition}[theorem]{Definition}
\theoremstyle{definition}

\newcommand{\beq}{\begin{equation}}
\newcommand{\eeq}{\end{equation}}
\newcommand{\beqs}{\begin{equation*}}
\newcommand{\eeqs}{\end{equation*}}
\newcommand{\beqn}{\begin{eqnarray}}
\newcommand{\eeqn}{\end{eqnarray}}
\newcommand{\beqns}{\begin{eqnarray*}}
\newcommand{\eeqns}{\end{eqnarray*}}
\newcommand{\bdoc}{\begin{document}}
\newcommand{\edoc}{\end{document}}
\newcommand{\be}{\begin{enumerate}}
\newcommand{\ee}{\end{enumerate}}
\newcommand{\bdescr}{\begin{description}}
\newcommand{\edescr}{\end{description}}
\newcommand{\ba}{\begin{array}}
\newcommand{\ea}{\end{array}}
\newcommand{\intR}{\int_{\mathbb R^N}}
\newcommand{\R}{\mathbb R}
\newcommand{\RN}{\mathbb{R}^N}
\newcommand{\B}{\mathbb B}
\newcommand{\C}{\mathcal C}
\renewcommand{\H}{\mathcal H}
\renewcommand{\L}{\mathbb L}
\newcommand{\parallelsum}{\mathbin{\!/\mkern-5mu/\!}}
\newcommand{\e}{\varepsilon}
\newcommand{\SD}{\Sigma_D}
 \renewcommand{\(}{\left(}
\renewcommand{\)}{\right)}
\renewcommand{\[}{\left[}
\renewcommand{\]}{\right]}
\renewcommand{\appendixpagename}{\centering Appendix}

\begin{document}

\makeatletter
\def\author@andify{%
  \nxandlist {\unskip ,\penalty-1 \space\ignorespaces}%
    {\unskip {} \@@and~}%
    {\unskip , \penalty-2}%
}
\makeatother

\title[Overdertermined problems and relative Cheeger sets]{Overdetermined problems and relative Cheeger sets in unbounded domains}


\author{D\MakeLowercase{anilo} G\MakeLowercase{regorin} Afonso}

\address[Danilo Gregorin Afonso]{Dipartimento di Matematica Guido Castelnuovo, Sapienza Universit\`a di Roma, Piazzale Aldo Moro 5, 00185 Roma, Italy}
\email{gregorinafonso@mat.uniroma1.it}


\author{A\MakeLowercase{lessandro} Iacopetti}
\address[Alessandro Iacopetti]{Dipartimento di Matematica ``G. Peano", Universit\`a di Torino, Via Carlo Alberto 10, 10123 Torino, Italy}
\email{alessandro.iacopetti@unito.it}

\author{F\MakeLowercase{ilomena} Pacella}
\address[Filomena Pacella]{Dipartimento di Matematica Guido Castelnuovo, Sapienza Universit\`a di Roma, Piazzale Aldo Moro 5, 00185 Roma, Italy}
\email{pacella@mat.uniroma1.it}


\subjclass[2020]{35J25, 35N25, 49Q10, 53A10}

\keywords{Overdetermined elliptic problem, shape optimization in unbounded domains, cartesian graphs, constant mean curvature surfaces}

\thanks{\emph{Acknowledgements.} Research partially supported by Gruppo Nazionale per l'Analisi Matematica, la Pro\-ba\-bi\-li\-t\`a e le loro Applicazioni (GNAMPA) of the Istituto Nazionale di Alta Matematica (INdAM)}

\maketitle

\begin{abstract}
In this paper we study a partially overdetermined mixed boundary value problem for domains $\Omega$ contained in an unbounded set $\mathcal C$. We introduce the notion of Cheeger set relative to $\mathcal C$ and show that if a domain $\Omega \subset \mathcal C$ admits a solution of the overdetermined problem, then it coincides with its relative Cheeger set. We also study the related problem of characterizing constant mean curvature surfaces $\Gamma$ inside $\mathcal C$. In the case when $\mathcal C$ is a cylinder we obtain further results whenever the relative boundary of $\Omega$ or the surface $\Gamma$ is a graph on the base of the cylinder.
\end{abstract}

\vskip2em

\section{Introduction}
\label{sec:introduction}

Given an unbounded Lipschitz domain $\mathcal C \subset \mathbb R^N$, $N \geq 2$, we consider a bounded domain $\Omega \subset \mathcal C$ and denote   
by $\Gamma_\Omega$ its \textit{relative (to $\mathcal C$) boundary}, i.e., $\Gamma_\Omega = \partial \Omega \cap \mathcal C$. Then we set $\Gamma_{1, \Omega} = \partial \Omega \cap \partial \mathcal C$ and assume that $\mathcal H_{N - 1}(\Gamma_{1, \Omega}) > 0$, where $\mathcal H_{N - 1}$ denotes the $(N-1)$-dimensional Hausdorff measure.

We study the partially overdetermined mixed boundary value problem
\begin{equation}
\label{eq:overdetermined_problem}
\left\{
\begin{array}{rcll}
- \Delta u & = & 1 & \quad \text{ in } \Omega \\
u &= & 0 & \quad \text{ on } \Gamma_\Omega \\
\displaystyle \frac{\partial u}{\partial \nu} & = & 0 & \quad \text{ on } \Gamma_{1, \Omega} \\
\displaystyle \frac{\partial u}{\partial \nu} & = & - c < 0 & \quad \text{ on } \Gamma_\Omega
\end{array}
\right.
\end{equation}
where $c > 0$ is a constant and $\nu$ denotes the exterior unit normal, which is defined on the regular part of $\partial \Omega$. The question we address is to determine the domains $\Omega \subset \mathcal C$ for which \eqref{eq:overdetermined_problem} admits a solutions. Obviously, this depends on the given container $\mathcal C$.

A related question is to characterize the surfaces $\Gamma$ in the set $\mathcal C$ with constant mean curvature (CMC, in short). More precisely, we consider CMC smooth, bounded, connected, orientable $(N-1)$-dimensional manifolds $\Gamma$ inside $\mathcal C$ whose boundary is contained in $\partial \mathcal C$ and intersects $\partial \mathcal C$ orthogonally. The characterization of such surfaces also depends on the set $\mathcal C$.

Both questions were addressed in \cite{PacellaTralli2020} (see also \cite{PacellaTralli2021}) in the case when $\mathcal C$ is a cone. The overdetermined problem for more general equations in cones was considered in \cite{CiraoloRoncoroni2020}.

When the cone is convex, the results of \cite{PacellaTralli2020} completely characterize the domains $\Omega$ for which \eqref{eq:overdetermined_problem} admits a solutions, as well as the CMC surfaces which intersect $\partial \mathcal C$ orthogonally: they are either spherical sectors centered at the vertex of the cone or half-balls lying on a flat portion of $\partial \mathcal C$. For nonconvex cones a similar characterization has been proved, but only for CMC radial graphs (\cite{PacellaTralli2020, PacellaTralli2021}). On the other side, nonradiality results have been proved for a class of nonconvex cones (\cite{IacopettiPacellaWeth2021}).

In this paper we study the aforementioned questions when $\mathcal C$ is a general unbounded smooth set. More precise results will be obtained when the unbounded set $\mathcal C$ is a cylinder in $\mathbb R^N$ spanned by a smooth bounded domain $\omega \subset \mathbb R^{N - 1}$. We denote it by $\mathcal C_\omega$, i.e.,
\begin{equation}
\label{eq:def_C_omega}
\mathcal C_\omega \coloneqq \omega \times (0, + \infty) = \{x = (x', x_N) \in \mathbb R^N \ : \ x' \in \omega, \ x_N \in (0, + \infty)\}.
\end{equation} 
In this case we obtain a result about CMC surfaces $\Gamma = \Gamma_\varphi$ inside $\mathcal C_\omega$ which are graphs of smooth functions $\varphi$ defined on $\overline \omega$. We prove that if $\Gamma_\varphi$ meets $\partial \mathcal C$ orthogonally, then the mean curvature of $\Gamma_\varphi$ is necessarily $0$ and $\Gamma_\varphi$ is the graph of a constant function (see Proposition \ref{prop:CMC_surfaces}).

Then we consider the Lipschitz domain 
$$
\Omega_\varphi \coloneqq \{(x', x_N) \in \mathbb R^N \ : \ x' \in \omega, \ x_N < \varphi(x')\} \subset \mathcal C_\omega
$$
and study the overdetermined problem \eqref{eq:overdetermined_problem} in $\Omega_\varphi$. In view of Proposition \ref{prop:CMC_surfaces} and the results in cones, it would be natural to conjecture that the domains $\Omega_\varphi$ for which \eqref{eq:overdetermined_problem} admits a solution are the ones corresponding to a function $\varphi \equiv h$ for some $h > 0$, so that $\Omega_\varphi$ is the bounded cylinder $\omega \times (0, h)$. In particular, $\Gamma_\varphi$ would have zero mean curvature. We are not able to prove this but we show some partial results.

In the case when $\mathcal C$ is a general Lipschitz unbounded domain we introduce the definition of relative (to $\mathcal C)$ Cheeger set for a domain $\Omega \subseteq \mathcal C$. It generalizes the classical one which goes back to \cite{Cheeger1970} (see also \cite{Leonardi2015, Parini2011}). Then we show that, if $\mathcal C$ is convex, any bounded domain $\Omega \subset \mathcal C$ for which \eqref{eq:overdetermined_problem} admits a solution coincides with its relative Cheeger set. This is an interesting geometric property of $\Omega$ which sheds light on the connections between overdetermined and isoperimetric problems. We also show some properties of relative self-Cheeger sets.

We delay the precise statement of the results and the comments on the proofs to the corresponding sections.

This paper is organized as follows. In Section \ref{sec:CMC_cartesian_graphs} we study CMC surfaces in cylinders. In Section \ref{sec:relative_Cheeger_sets} we introduce the relative Cheeger problem and prove some properties of relative Cheeger sets. We also provide a bound on the first eigenvalue of the mixed boundary value problem for the Laplacian in terms of the Cheeger constant. The overdetermined problem is studied in Section \ref{sec:the_overdetermined_problem}.


\vskip2em

\section{CMC cartesian graphs}
\label{sec:CMC_cartesian_graphs}
Let $\omega\subset\R^{N-1}$ be a smooth bounded domain and let $\mathcal{C}_\omega\subset \R^N$ be the half-cylinder spanned by $\omega$, i.e., $\mathcal{C}_\omega=\omega\times(0,+\infty)$. We denote by $\partial\C_\omega^+$ the lateral part of $\partial\C_\omega$, i.e., $$\partial\C_\omega^+ \coloneqq \partial\C_\omega\cap\{x_N>0\}.$$
Let $\varphi\in C^2(\overline\omega)$ be such that $\varphi>0$ in $\overline\omega$ and consider the associated cartesian graph 
$$\Gamma_{\varphi}\coloneqq \{(x^\prime,x_N)\in \R^N \ : \ x' \in \omega, \ x_N=\varphi(x')\}.$$
By construction we have that $\Gamma_{\varphi} \subset \mathcal{C}_\omega$ and $\Gamma_\varphi$ meets $\partial\C_\omega$ only on $\partial\C_\omega^+$. The main result of this section is the following:

\begin{proposition}
\label{prop:CMC_surfaces}
Assume that $\Gamma_\varphi$ is a CMC surface which meets $\partial\mathcal{C}_\omega^+$ orthogonally. Then $\Gamma_\varphi$ is a minimal hypersurface and $\varphi$ is a constant function.
\end{proposition} 

\begin{proof}
Since $\Gamma_\varphi$ is a smooth cartesian graph, then the mean curvature of $\Gamma_\varphi$ (with respect to the exterior unit normal $\nu_{\Gamma_\varphi}$) is given by
\begin{equation}\label{eq1sect2}
H_{\Gamma_\varphi}=-\frac{1}{N-1}\mathrm{div}\left(\frac{\nabla \varphi}{\sqrt{1+|\nabla \varphi|^2}}\right)
\end{equation}
and, as $H_{\Gamma_\varphi}$ is constant, then integrating \eqref{eq1sect2} over $\omega$ we deduce that
\begin{equation}\label{eq1bissect2}
(N-1) H_{\Gamma_\varphi} |\omega|=  -\int_\omega\mathrm{div}\left(\frac{\nabla \varphi}{\sqrt{1+|\nabla \varphi|^2}}\right) \ dx^\prime.
\end{equation}
On the other hand, by the divergence theorem we have
\begin{equation}\label{eq2sect2}
\int_\omega\mathrm{div}\left(\frac{\nabla \varphi}{\sqrt{1+|\nabla \varphi|^2}}\right) \ dx^\prime = \int_{\partial\omega} \frac{\nabla \varphi \cdot \nu_{\partial\omega}}{\sqrt{1+|\nabla \varphi|^2}} \ d\sigma,
\end{equation}
where $\nu_{\partial\omega}$ is the unit outward normal to $\partial\omega$. Let us observe that since $\Gamma_\varphi$ meets $\partial\mathcal{C}_\omega^+$ orthogonally, then for all $x\in\overline\Gamma_\varphi\cap \partial\mathcal{C}_\omega$ it holds
\begin{equation}\label{eq3sect2}
\nu_{\Gamma_\varphi}(x)\cdot \nu_{\partial\mathcal{C}_\omega^+}(x) =0.
\end{equation}
Moreover, as $\Gamma_\varphi$ is a cartesian graph,  for all $x=(x^\prime,x_N)\in\overline\Gamma_\varphi$ we have
\begin{equation}\label{eq4sect2} 
\nu_{\Gamma_\varphi}(x)=\nu_{\Gamma_\varphi}(x', \varphi(x'))=\frac{1}{\sqrt{1+|\nabla \varphi(x^\prime)|^2}}\left[\begin{array}{c}
-\nabla\varphi(x^\prime)\\
1\end{array}\right].
\end{equation}
In particular, since $\partial\mathcal{C}_\omega^+$ is a cylinder spanned by $\partial\omega$ we have $\nu_{\partial\mathcal{C}_\omega^+}(x)=(\nu_{\partial\omega}(x^\prime),0)$ for all $x=(x^\prime,x_N) \in  \partial\mathcal{C}_\omega^+$. From these considerations, \eqref{eq3sect2} and \eqref{eq4sect2} we readily obtain that 
\begin{equation}\label{eq5sect2} 
\nabla \varphi(x^\prime) \cdot \nu_{\partial\omega}(x^\prime)=0 \ \ \ \forall x\in\partial\omega.
\end{equation}
Hence, combining \eqref{eq1bissect2}, \eqref{eq2sect2} and \eqref{eq5sect2} we deduce that
$$H_{\Gamma_\varphi}=0.$$
To conclude it remains to prove that $\varphi$ is a constant function. Since $H_{\Gamma_\varphi} = 0$, then from \eqref{eq1sect2}, integrating by parts and taking into account \eqref{eq5sect2} it follows that
$$ \int_\omega \frac{\nabla \varphi \cdot \nabla \psi}{{\sqrt{1+|\nabla \varphi|^2}}} \ dx^\prime = \int_{\partial \omega} \psi  \frac{\nabla \varphi \cdot\nu_{\partial\omega}}{{\sqrt{1+|\nabla \varphi|^2}}} \ d\sigma=0 \ \ \ \forall \psi \in C^1(\overline\omega).$$
Finally, choosing $\psi=\varphi$ we readily obtain that
$$ \int_\omega \frac{|\nabla \varphi|^2}{{\sqrt{1+|\nabla \varphi|^2}}} \ dx^\prime = 0,$$
which implies $|\nabla \varphi|^2 =0$ on $\omega$. As $\omega$ is connected, we get that $\varphi$ is a constant function.
\end{proof}

Finally we prove the following Minkowski formula for graphs, which is interesting in itself and will be used in Section \ref{sec:the_overdetermined_problem}.

\begin{proposition}
\label{prop:Minkoski_formula_for_graphs}
Assume that $\Gamma_\varphi$ meets $\partial\mathcal{C}_\omega^+$ orthogonally. Then
\begin{equation}\label{eqthesisMink}
\int_{\Gamma_\varphi} H_{\Gamma_\varphi} \langle x_N e_N, \nu\rangle \ d\sigma= \frac{1}{N-1} \int_\omega \frac{|\nabla\varphi|^2}{\sqrt{1+|\nabla\varphi|^2}} \ dx',
\end{equation}
where $e_N = (0,\ldots,0,1)\in \R^N$.
\end{proposition}

\begin{proof}
Since $\Gamma_\varphi$ is the cartesian graph associated to $\varphi$ we have 
\begin{equation}\label{eq:dsigma}
d\sigma=\sqrt{1+|\nabla \varphi|^2} \ dx^\prime. 
\end{equation}

Exploiting \eqref{eq1sect2} and \eqref{eq4sect2} and integrating by parts we get that
\begin{align}
\int_{\Gamma_\varphi} H_{\Gamma_\varphi} \langle x_N e_N, \nu\rangle \ d\sigma 
& = -\frac{1}{N-1} \int_\omega \mathrm{div}\left(\frac{\nabla \varphi}{\sqrt{1+|\nabla \varphi|^2}}\right)  \varphi \ dx^\prime \nonumber \\ 
& = \frac{1}{N-1} \left(\int_\omega\frac{|\nabla \varphi|^2}{\sqrt{1+|\nabla \varphi|^2}} \ dx^\prime - \int_{\partial \omega} \varphi  \frac{\nabla \varphi \cdot\nu_{\partial\omega}}{\sqrt{1+|\nabla \varphi|^2}} \ d\sigma_{\partial \omega}\right) \label{eq1:Mink}
\end{align}
Finally, since $\Gamma_\varphi$ intersects $\partial\C_\omega^+$ orthogonally, then from \eqref{eq5sect2} we infer that the last integral in \eqref{eq1:Mink} is zero, and \eqref{eqthesisMink} readily follows.
\end{proof}

\vskip2em

\section{Relative Cheeger sets}
\label{sec:relative_Cheeger_sets}

Let $\mathcal C \subset \mathbb R^N$ be an unbounded Lipschitz domain. Let us recall the definition of relative perimeter.

\begin{definition}
The \textit{relative (to $\mathcal C$) perimeter} of a set $E \subset \mathcal C$ is
$$
P_{\mathcal C}(E) \coloneqq |D \chi_E|_{\mathcal C}(E) = \sup \left\{\int_E \divergence \psi \ dx \ : \ \psi \in C_c^1(\mathcal C, \mathbb R^N), \ \|\psi\|_\infty \leq 1 \right\},
$$
where $\chi_E$ is the characteristic function of $E$. Of course, a set is said to be of \textit{finite (relative) perimeter} when  $P_{\mathcal C}(E) < + \infty$.
\end{definition}

Note that if $\Gamma_E:=\partial E\cap \mathcal C$ is Lipschitz, then
\begin{equation}
\label{eq:relative_perimeter_of_Omega}
P_{\mathcal C} (E) = \mathcal H_{N - 1}(\Gamma_E) = \int_{\Gamma_E} 1 \ d\sigma.
\end{equation}

Let $|E|$ denote the Lebesgue measure of $E$, which will also be called the volume of $E$. Let $\Omega \subseteq \mathcal C$ and set, as in Section 1, $\Gamma_\Omega \coloneqq \partial \Omega \cap \mathcal C$.

\begin{definition}
The \textit{relative (to $\mathcal C$) Cheeger constant} of $\Omega \subseteq \mathcal C$ is
$$
h_{\mathcal C}(\Omega) \coloneqq \inf_{E \subseteq \Omega} \frac{P_{\mathcal C}(E)}{|E|}.
$$
The sets which attain the minimum will be called the \textit{relative (to $\mathcal C$) Cheeger sets} of $\Omega$. If $\Omega$ itself is a Cheeeger set, then $\Omega$ is said to be \textit{self-Cheeger}.
\end{definition}

The problem of finding $h_{\mathcal C}(\Omega)$ and the associated Cheeger set (or sets) has many interesting motivations and applications. The Cheeger constant $h(\Omega)$ first appeared in a bound for the first eigenvalue of the Laplacian on manifolds (\cite{Cheeger1970}); see \cite{KawohlFridman2003} for the link with the spectral theory for the $p$-Laplacian. See \cite{Keller1980} and the references therein for fracture problems in mechanics of materials; \cite{IonescuLachandRobert2005} for a landslide problem; \cite{CasellesFaccioloMeinhardt2009} for a generalization of the Cheeger problem in image processing. Regarding the Cheeger problem on its own, see \cite{BuenoErcole2011} for a study via $p$-torsion functions; \cite{KawohlLachandRobert2006} for a characterization of Cheeger sets of convex plane domains. The surveys \cite{Parini2011, Leonardi2015} provide a nice overview.

We now show some properties of relative Cheeger sets.

\begin{proposition}[Existence]
Let $\Omega \subset \mathcal C$ be a bounded domain such that $\Gamma_\Omega$ is Lipschitz. Then there exists at least one Cheeger set for $\Omega$.
\end{proposition}

\begin{proof}
The proof is a straightforward adaptation of the proof of Proposition 3.1 in \cite{Parini2011}. The important thing to note is that the test functions $\varphi$ have compact support in $\mathcal C$, so they don't ``see" the boundary of the container.
\end{proof}

\begin{proposition}
Let $\mathcal C$ be a cone with vertex at the origin, i.e., $\mathcal C = \{tx \ : \ x \in D, \ t \in (0, + \infty)\}$ where $D$ is a domain on the unit sphere $S^{N - 1}$. Let $\Omega \subset \mathcal C$ be a bounded domain such that $\Gamma_\Omega$ is Lipschitz. If $E$ is a Cheeger set for $\Omega$, then $\partial E \cap \Gamma_\Omega \neq \emptyset$. 
\end{proposition}

\begin{proof}
We adapt the proof of \cite[Proposition 3.5]{Parini2011}. Note that the fact that the cone is invariant by dilation  plays a crucial role in the proof.

For the sake of contradiction, suppose that $d(\partial E, \Gamma_\Omega) \geq \delta$ for some $\delta > 0$. Then we can find some $t > 1$ such that the set
$$
tE\coloneqq \{x \in \mathcal C \ : \ t^{-1}x \in E\}
$$
is still in $\Omega$. By the change of variables formula we have
$$
|tE| = \int_{tE} 1 \ dx = \int_E t^N \ dx = t^N |E|.
$$
To compute the perimeter one proceeds as follows. By the definition of relative perimeter, the chain rule and the change of variables formula we get 
\begin{align*}
P_{\mathcal C}(tE)
& = \sup\left\{\int_{tE} \divergence \psi (x) \ dx \ : \ \psi \in C_c^1(\mathcal C, \mathbb R^N), \ \|\psi\|_\infty \leq 1 \right\} \\
& = \sup \left\{\int_E \divergence (\psi(t^{-1}x)) t^N \ dx \ : \ \psi \in C_c^1(\mathcal C, \mathbb R^N), \ \|\psi\|_\infty \leq 1 \right\} \\
& = t^{N - 1} P_{\mathcal C}(E).
\end{align*} 
Observe that the set of functions where we compute the supremum is the same, but when computing the divergence we must take into account the dilation of the argument.

Now it remains to observe that
$$
\frac{P_{\mathcal C}(tE)}{|tE|} = \frac{t^{N-1}}{t^N} \frac{P_{\mathcal C}(E)}{|E|} < h_\mathcal C(\Omega),
$$
which contradicts the definitions of (relative) Cheeger constant and Cheeger set.
\end{proof}

We can prove a similar result in cylinders.

\begin{proposition}
Let $\omega \subset \mathbb R^{N - 1}$ and let $\mathcal C_\omega$ be the cylinder spanned by $\omega$. Assume that $\Omega \subset \mathcal C_\omega$ is a bounded domain such that $\Gamma_\Omega$ is a connected surface whose projection into $\mathbb R^{N - 1}$ is exactly $\omega$. If some relative Cheeger set $E_f$ of $\Omega$ is defined by the graph of a function $f$ on $\omega$, i.e., $E_f = \{(x', x_N) \in \mathcal C_\omega \ : \ 0 < x_N < f(x')\}$, then $\Gamma_{E_f} \cap \Gamma_\Omega \neq \emptyset$.
\end{proposition}

\begin{proof}
Note that $\Gamma_{E_f} = \partial E_f \cap \mathcal C$ is just the graph $\Gamma_f$ of the function $f$. If $\Gamma_f \cap \Gamma_\Omega = \emptyset$, considering the domain $E_{f + \delta}$ defined by the function $f + \delta$, for $\delta > 0$ sufficiently small so that $E_{f + \delta } \subseteq \Omega$, we have  
$$
\frac{\mathcal H_{N - 1}(\Gamma_{f + \delta})}{|E_{f + \delta}|} = \frac{\mathcal H_{N - 1}(\Gamma_f)}{|E_{f + \delta}|} < \frac{\mathcal H_{N - 1}(\Gamma_f)}{|E_f|},
$$
contradicting the definition of a Cheeger set.
\end{proof}

Note that in the previous proof we can allow $f$ to be zero at some strict subset of $\omega$.

\begin{remark}\label{rem:smoothness}
The Cheeger constant can be obtained by minimization in the class of smooth subdomains of $\Omega$, see the proof of Theorem 1 in \cite{LionsPacella1990} and the references therein.
\end{remark}

When $\mathcal C$ is a cylinder and $\Omega_\varphi$ is the domain defined by the graph of a function $\varphi \in C^2(\overline \omega)$ we have the following result:

\begin{theorem}
\label{thm:graph_is_Cheeger_iff_is_constant}
Let $\omega \subset \mathbb R^{N - 1}$ be a smooth bounded domain and consider the cylinder $C_\omega$ spanned by $\omega$. Let $\varphi \in C^2(\overline \omega)$ be a positive function and assume that $\Gamma_\varphi$ meets $\partial \mathcal C_\omega$ orthogonally. If $\Omega_\varphi$ is self-Cheeger, then
\begin{equation}
\label{eq:bound_H_of_self_Cheeger}
H_{\Gamma_\varphi} \leq \frac{P_{\mathcal C_\omega}(\Omega_\varphi)}{(N - 1)|\Omega_\varphi|}
\end{equation}
\end{theorem}

\begin{proof}
For any $C^1$ function $v \leq 0$ on $ \overline \omega$, the relative perimeter of the perturbed domain $\Omega_{\varphi + tv}$ is given by (recall \eqref{eq:relative_perimeter_of_Omega} and \eqref{eq:dsigma}):
\begin{align}
p_v(t) 
& \coloneqq P_{\mathcal C}(\Omega_{\varphi + tv}) = \mathcal H_{N - 1}(\Gamma_{\varphi + tv}) \nonumber \\
& = \int_\omega \sqrt{1 + |\nabla \varphi + t \nabla v|^2} \ dx' \label{eq:p_w(t)}.
\end{align}
We also define 
\begin{equation}
\label{eq:V(t)}
V_v(t) \coloneqq |\Omega_{\varphi + tw}| = \int_\omega \varphi + t v \ dx'.
\end{equation}
Here we are considering $t$ small enough for these quantities to make sense in the container, that is, $\varphi + tv \geq 0$.

If $\Omega_\varphi$ is self-Cheeger, then
\begin{equation}
\label{eq:perturbation_derivative_of_self_Cheeger}
\left. \frac{d}{dt} \left( \frac{p_v(t)}{V_v(t)}\right)\right|_{t = 0} = \left. \frac{p_v'(t)V_v(t) - p_v(t) V_v'(t)}{(V_v(t))^2} \right|_{t = 0} \geq 0
\end{equation}
for all negative $v \in C^1(\overline \omega)$.

Since $(V_v(t))^2$ is always positive, we have that
\begin{equation}
\label{eq:perturbation_derivative_of_self_Cheeger_2}
p_v'(0) V_v(0) - p_v(0) V_v'(0) \geq 0.
\end{equation}

Simple computations yield
\begin{equation}
\label{eq:p_v'(t)}
p_v'(t) = \int_\omega \frac{ \nabla \varphi \cdot \nabla v + t |\nabla v|^2}{\sqrt{1 + |\nabla \varphi + t \nabla v|^2}} \ dx'
\end{equation}
and
\begin{equation}
\label{eq:V_v'(t)}
V_v'(t) = \int_\omega v \ dx'.
\end{equation}
Then
\begin{align}
\label{eq:p_v'(0) V_v(0) - p_v(0) V_v'(0)}
p_v'(0) V_v(0) - p_v(0) V_v'(0)
& = \int_\omega \frac{\nabla \varphi \cdot \nabla v}{\sqrt{1 + |\nabla \varphi|^2}} \ dx' |\Omega_\varphi| - P_{\mathcal C_\omega}(\Omega_\varphi) \int_\omega v \ dx'.
\end{align}

Integrating by parts and using the divergence theorem we obtain
\begin{align}
\int_\omega \frac{\nabla \varphi \cdot \nabla v}{\sqrt{1 + |\nabla \varphi|^2}} \ dx' 
& = \int_\omega \divergence \left(v \frac{\nabla \varphi}{\sqrt{1 + |\nabla \varphi|^2}} \right) \ dx' - \int_\omega v \divergence \left(\frac{\nabla \varphi}{\sqrt{1 + |\nabla \varphi|^2}} \right) \ dx' \nonumber \\
& = \int_{\partial \omega} v \frac{\nabla \varphi \cdot \nu_{\partial \omega}}{\sqrt{1 + |\nabla \varphi|^2}} \ dx' - \int_\omega v \divergence \left(\frac{\nabla \varphi}{\sqrt{1 + |\nabla \varphi|^2}} \right) \ dx' \nonumber \\ 
& = - \int_\omega v \divergence \left(\frac{\nabla \varphi}{\sqrt{1 + |\nabla \varphi|^2}} \right) \ dx, \label{eq:getting_the_curvature}
\end{align}
since $\Gamma_\varphi$ meets $\partial \mathcal C_\omega$ orthogonally (see the proof of Proposition \ref{prop:CMC_surfaces}). Then, substituting \eqref{eq:getting_the_curvature} into \eqref{eq:p_v'(0) V_v(0) - p_v(0) V_v'(0)} and taking into account \eqref{eq:perturbation_derivative_of_self_Cheeger_2} we get
\begin{equation}
\int_\omega v \left(|\Omega_\varphi|\divergence\left(\frac{\nabla \varphi}{\sqrt{1 + |\nabla \varphi|^2}}\right) + P_{\mathcal C_\omega}(\Omega_\varphi)\right) \ dx' \leq 0
\end{equation}
for every negative $v \in C^1(\overline\omega)$. Hence
\begin{align}
0 \leq |\Omega_\varphi|\divergence\left(\frac{\nabla \varphi}{\sqrt{1 + |\nabla \varphi|^2}}\right) + P_{\mathcal C_\omega}(\Omega_\varphi) 
& = - (N - 1) H_{\Gamma_\Omega} |\Omega_\varphi| + P_{\mathcal C_\omega}(\Omega_\varphi),
\end{align}
from which \eqref{eq:bound_H_of_self_Cheeger} readily follows.
\end{proof}

\begin{theorem}
\label{thm:if_nabla_bounded_below_then_not_Cheeger}
Let $\omega \subset \mathbb R^{N - 1}$ be a smooth bounded domain and consider the cylinder $C_\omega$ spanned by $\omega$. Let $\varphi \in C^2(\overline \omega)$ be a positive function. Assume that $|\nabla \varphi| \geq \delta$ on $\overline\omega$, for some $\delta > 0$. Then $\Omega_\varphi$ cannot be self-Cheeger.
\end{theorem}

\begin{proof}
Again we apply the idea of perturbing $\varphi$, but now the aim is to find a specific function $v \leq 0$ in $C^1(\overline \omega)$ for which
\begin{equation}
\label{eq:perturbation_derivative}
\left. \frac{d}{dt} \left( \frac{p_v(t)}{V_v(t)}\right)\right|_{t = 0} = \left. \frac{p_v'(t)V_v(t) - p_v(t) V_v'(t)}{(V_v(t))^2} \right|_{t = 0} < 0.
\end{equation}
Indeed this will imply that
\begin{equation}
\frac{P_{\mathcal C}(\Omega_{\varphi + tv})}{|\Omega_{\varphi + tv}|} < \frac{P_{\mathcal C}(\Omega)}{|\Omega|}
\end{equation}
for every $t < \varepsilon$ for some $\varepsilon > 0$.

As $(V_v(t))^2$ is always positive, the problem reduces to finding $v$ such that
$$
p_v'(0) V_v(0) - p_v(0) V_v'(0) < 0.
$$

Take
$$
v(x') = - e^{\alpha \varphi(x')}
$$
where $\alpha > 0$ is a constant to be chosen later.
Then
$$
\nabla v(x') = - \alpha e^{\alpha \varphi(x')} \nabla \varphi(x')
$$
and
\begin{align}
p_v'(0) V_v(0) - p_v(0) V_v'(0)
& = \int_\omega - \alpha |\Omega_\varphi| e^{\alpha \varphi} \frac{|\nabla \varphi|^2}{\sqrt{1 + |\nabla \varphi|^2}} \ dx' + \int_\omega P_{\mathcal C}(\Omega_\varphi) e^{\alpha \varphi} \ dx'\nonumber \\
& \leq \int_\omega e^{\alpha \varphi} \left(P_{\mathcal C}(\Omega_\varphi) - \alpha |\Omega_\varphi| \frac{|\nabla \varphi|^2}{\sqrt{1 + |\nabla \varphi|^2}} \right) \ dx' \label{eq:nabla_varphi_bounded_then_not_self_Cheeger}.\nonumber
\end{align}
Since $|\nabla \varphi|$ is bounded on $\overline\omega$ (because $\varphi$ is smooth) and $|\nabla \varphi|\geq\delta$ on $\overline\omega$, choosing $\alpha$ big enough such that
$$
P_{\mathcal C}(\Omega_\varphi) - \alpha |\Omega_\varphi| \frac{\delta^2}{\sqrt{1 + \|\nabla \varphi\|_\infty^2}} < 0
$$
we easily conclude.
\end{proof}

We now go back to the case when $\mathcal C$ is a general unbounded Lipschitz domain and conclude this section with an application.

Assume that $\Gamma_\Omega$ is smooth and consider the eigenvalue problem with mixed boundary conditions:
\begin{equation}
\label{eq:mixed_eigenvalue_problem}
\left\{
\begin{array}{rcll}
- \Delta u & = & \lambda u & \quad \text{ in } \Omega \\
u & = & 0 & \quad \text{ on } \Gamma_\Omega \\
\displaystyle \frac{\partial u}{\partial \nu} & = & 0 & \quad \text{ on } \Gamma_{1, \Omega}
\end{array}
\right.
.
\end{equation}

Let $H_0^1(\Omega \cup \Gamma_{1, \Omega})$ be the closure in $H^1(\Omega)$ of the space $C_c^1(\Omega \cup \Gamma_{1, \Omega})$. It is the natural space to study problem \eqref{eq:mixed_eigenvalue_problem}, for it coincides with the space of functions in $H^1(\Omega)$ whose trace vanishes on $\Gamma_\Omega$.

The spectral theory for this problem is analogous to the one for the Dirichlet-Laplacian\footnote{See \cite[Section 1.4]{DamascelliPacella2019} for the spectral theory of mixed boundary value problems in a more general setting.}. In particular, the eigenvalues are all positive, form an increasing divergent  sequence and the first eigenfunction is positive. Moreover, the first eigenvalue is characterized by
$$
\lambda_1(\Omega) = \min_{v \in H_0^1(\Omega \cup \Gamma_{1, \Omega}) \setminus \{0\}} \frac{\displaystyle \int_\Omega |\nabla v|^2 \ dx}{\displaystyle \int_\Omega v^2 \ dx} .
$$


\begin{theorem}
Let $\Omega \subset \mathcal C$ be a bounded domain such that $\Gamma_\Omega$ is Lipschitz. Then
$$
\lambda_1(\Omega) \geq \frac{h_{\mathcal C}^2(\Omega)}{4}.
$$
\end{theorem}

\begin{proof}
We adapt the proof presented for the Dirichlet-$p$-Laplacian in the Appendix of \cite{LeftonWei1997}

Let $w \in C_c^\infty(\Omega \cup \Gamma_{1, \Omega})$ and $E_t = \{x \in \Omega \ : \ w(x) > t\}$. By the classical coarea formula, by the definition of $h_{\mathcal C}(\Omega)$ and Cavalieri's principle, we have
\begin{align*}
\int_\Omega |\nabla w| \ dx
& = \int_{- \infty}^\infty \mathcal H_{N - 1}(w^{-1}(t)) \ dt \\
& = \int_{- \infty}^\infty \frac{\mathcal H_{N - 1}(w^{-1}(t))}{|E_t|} |E_t| \ dt \\
& \geq h_{\mathcal C}(\Omega) \int_\Omega |w| \ dx.
\end{align*}
Hence
\begin{equation}
\label{eq:ineq_w_C_c}
h_{\mathcal C}(\Omega) \leq \frac{\displaystyle \int_\Omega |\nabla w| \ dx}{\displaystyle \int_\Omega w \ dx} \quad \forall w \in C_c^\infty(\Omega \cup \Gamma_{1, \Omega}),
\end{equation}
and by density \eqref{eq:ineq_w_C_c} holds also in $W_0^{1, 1}(\Omega \cup \Gamma_{1, \Omega})$.

Now let $v \in H_0^1(\Omega \cup \Gamma_{1, \Omega})$. By H\"{o}lder's inequality, it follows that $v^2 \in W_0^{1, 1}(\Omega \cup \Gamma_{1, \Omega})$. Indeed
$$
\int_\Omega |\nabla (v^2)| \ dx = 2 \int_\Omega |v| |\nabla v| \ dx \leq 2 \|v\|_2 \|\nabla v\|_2.
$$
By \eqref{eq:ineq_w_C_c} it follows that
$$
h_{\mathcal C}(\Omega) \leq \frac{\displaystyle \int_\Omega |\nabla(v^2)| \ dx}{\displaystyle \int_\Omega v^2 \ dx} \leq 2 \frac{\|v\|_2 \|\nabla v\|_2}{\|v\|_2^2} = 2 \frac{\|\nabla v\|_2}{\|v\|_2}.
$$
We conclude the proof by taking into account the variational characterization of the first eigenvalue $\lambda_1(\Omega)$.
\end{proof}

\vskip2em

\section{The overdetermined problem}
\label{sec:the_overdetermined_problem}

Let $\mathcal C$ be an unbounded Lipschitz domain and $\Omega \subset \mathcal C$ be a bounded domain with smooth relative boundary $\Gamma_\Omega$.

%

For a solution $u$ of \eqref{eq:overdetermined_problem} we define
\begin{equation}
\label{eq:P_function}
P(x) = |Du(x)|^2 + \frac{2}{N} u(x), \quad x \in \overline \Omega.
\end{equation}
This function (sometimes called P-function) is often used in the study of overdetermined problems (\cite{Weinberger1971, PacellaTralli2020, FragalaGazzolaKawohl2006}).

We can give a bound on the curvature of $\Gamma_\Omega$, with the aid of the following

\begin{lemma}
\label{lem:P_function}
If $\mathcal C$ is convex and $u$ is a solution of \eqref{eq:overdetermined_problem}, then either $P \equiv c^2$ in $\overline \Omega$ or $\frac{\partial P}{\partial \nu} > 0$ on $\Gamma_\Omega$.
\end{lemma}

\begin{proof}
Direct computations yield
\begin{equation}
\label{eq:DP}
DP = 2 D^2u Du + \frac{2}{N} Du
\end{equation}
and
\begin{equation}
\label{eq:DeltaP}
\Delta P = 2 \|D^2u\|^2 + 2 \langle Du, D(\Delta u) \rangle - \frac{2}{N} = 2 \left(\|D^2u\|^2 - \frac{(\Delta u)^2}{N} \right),
\end{equation}
since $\Delta u = -1$. Here $D^2u$ denotes the Hessian matrix of $u$ and $\|D^2u\|^2$ is the sum of the squares of the elements of $D^2u$. By \cite[inequality (2.4)]{PacellaTralli2020} it then follows that
\begin{equation}
\label{eq:DeltaP_geq_0}
\Delta P \geq 0 \ \ \hbox{in $\Omega$.}
\end{equation}
Moreover, by the boundary conditions for $u$ we obtain the following boundary conditions for $P$:
\begin{equation}
\label{eq:boundary_conditions_for_P}
P \equiv c^2 \ \text{ on } \Gamma_\Omega, \qquad \frac{\partial P}{\partial \nu} = 2 \langle D^2u Du, \nu \rangle \ \text{ on } \Gamma_{1, \Omega}.
\end{equation}

From the convexity assumption on $\mathcal C$ we have that the second fundamental form $h(\cdot, \cdot)$ on $\partial \mathcal C$ is positive semidefinite at any regular point. On the other hand, by the Neumann condition $\frac{\partial u}{\partial \nu} = 0$ on $\Gamma_{1, \Omega}$ we deduce that $Du$ is tangent to $\Gamma_{1, \Omega}$. Thus, denoting by $N_u$ the vector field differentiating along the direction $Du$, i.e., 
$$
N_u = \sum_{k = 1}^N u_k(x) \frac{\partial}{\partial x_k},
$$
we obtain
$$
0 = N_u(\langle Du, \nu \rangle) = \sum_{j, k = 1}^N u_k u_{kj} \nu_j + h(Du, Du) \geq \langle D^2u Du, \nu \rangle.
$$
Hence the function $P$ satisfies
\begin{equation}
\label{eq:differential_inequality_for_P}
\left\{
\begin{array}{rcll}
\Delta P & \geq & 0 & \quad \text{ in } \Omega \\
P & = & c^2 & \quad \text{ on } \Gamma_\Omega \\
\displaystyle \frac{\partial P}{\partial \nu} & \leq & 0 & \quad \text{ on } \Gamma_{1, \Omega}
\end{array}
.
\right.
\end{equation}

By the maximum principle for the mixed boundary value problem (see \cite[Section 1.2.1]{DamascelliPacella2019} or \cite[Corollary 2.3]{PacellaTralli2020}) we obtain $P \leq c^2$ in $\Omega$. Then the strong maximum principle applies, so that either $P \equiv c^2$ in $\Omega$ or $P < c^2$ in $\Omega$. In this last case, by Hopf's Lemma we get $\frac{\partial P}{\partial \nu} > 0$ on $\Gamma_\Omega$.
\end{proof}

Let $H(x)$ denote the mean curvature at a point $x \in \Gamma_\Omega$.

\begin{proposition}
\label{prop:bound_on_H_in_Gamma}
If $\mathcal C$ is convex and there exists a solution for \eqref{eq:overdetermined_problem} in $\Omega$, then either
\begin{equation}
\label{eq:Bound_H_from_P_1}
H(x) < \frac{1}{Nc}
\end{equation}
or
\begin{equation}
\label{eq:const_H_from_P}
H(x) \equiv \frac{1}{Nc}
\end{equation}
on for every $x \in \Gamma_\Omega$.
\end{proposition}

\begin{proof}
The proof of \cite[Lemma 3.3]{FragalaGazzolaKawohl2006} can be easily adapted to our setting. Denoting for brevity $u_\nu=\frac{\partial u}{\partial \nu}$, exploiting that $- \Delta u = 1$ and arguing as in  \cite[Lemma 3.3]{FragalaGazzolaKawohl2006} we have
\begin{equation}
\label{eq:FGK_12}
u_{\nu \nu} - (N - 1)cH(x) = -1 \ \ \hbox{on $\Gamma_\Omega$}.
\end{equation}

Consider the two possible cases given by Lemma \ref{lem:P_function}. In the inequality case,
\begin{equation}
\label{eq:FGK_13}
P_\nu = 2 u_\nu u_{\nu \nu} + \frac{2}{N} u_\nu > 0.
\end{equation}
Dividing \eqref{eq:FGK_13} by $2 u_\nu$ and combining with \eqref{eq:FGK_12} we obtain 
$$
H(x) < \frac{1}{Nc}.
$$ 
The case of equality in Lemma \ref{lem:P_function} develops similarly.
\end{proof}

As anticipated, we show that solutions of the overdetermined problem are (relatively) self-Cheeger:

\begin{theorem}
\label{thm:Omega_is_Cheeger}
If $\mathcal C$ is convex and there exists a solution $u$ for \eqref{eq:overdetermined_problem} in $\Omega$, then $|Du| \leq c$ in $\overline \Omega$, $c$ as in \eqref{eq:overdetermined_problem}, and $\Omega$ is self-Cheeger.
\end{theorem}

\begin{proof}
Recall from the proof of Lemma \ref{lem:P_function} that
$$
P(x) = |Du(x)|^2 + \frac{2}{N} u(x) \leq c^2
$$
in $\overline \Omega$. By the maximum principle we have that $u$ is positive, and thus we obtain that $|Du| \leq c$ in $\overline \Omega$. It then follows that for any smooth subdomain $E \subset \Omega$ (recall Remark \ref{rem:smoothness}) we have
$$
|E| = \int_E - \Delta u \ dx = \int_{\partial E} \frac{\partial u}{\partial \nu} \ d\sigma \leq c \mathcal H_{N - 1}(\partial E \cap \mathcal C)
$$
because $\frac{\partial u}{\partial \nu} = 0$ on $\Gamma_{1, \Omega}$.

Now, integrating \eqref{eq:overdetermined_problem} in $\Omega$ we have
\begin{equation}
\label{eq:|Omega|_equal_c_|Gamma|}
|\Omega| = - \int_\Omega \Delta u \ dx = - \int_{\partial \Omega} \frac{\partial u}{\partial \nu} \ d \sigma = \int_{\Gamma_\Omega} c \ d \sigma = c \mathcal H_{N - 1}(\Gamma_\Omega).
\end{equation}
Hence
$$
\frac{\mathcal H_{N - 1}(\Gamma_\Omega)}{|\Omega|} = \frac{1}{c} \leq \frac{\mathcal H_{N - 1}(\partial E \cap \mathcal C)}{|E|}.
$$
The proof is complete.
\end{proof}

\begin{remark}
The previous theorem allows to study the geometrical properties of the bounded domains $\Omega$ which admit a solution of the overdetermined problem \eqref{eq:overdetermined_problem}, by looking at those of the self-Cheeger sets. In particular, we point out those provided by Theorems \ref{thm:graph_is_Cheeger_iff_is_constant} and \ref{thm:if_nabla_bounded_below_then_not_Cheeger}.
\end{remark}

We conclude with a result which applies when $\Omega=\Omega_\varphi$ is defined by the graph of a function $\varphi$ as in Section 3. It gives a bound for the gradient of $\varphi$, whenever $\Omega_\varphi$ admits a solution of the overdetermined problem \eqref{eq:overdetermined_problem}.

\begin{proposition}
Let $\omega \subset \mathbb R^{N - 1}$ be a convex smooth bounded domain and $\varphi \in C^2(\overline \omega)$ be a positive function. If $\Gamma_\varphi$ meets $\partial \mathcal C_\omega$ orthogonally and $\Omega_\varphi$ admits a solution to the overdetermined problem \eqref{eq:overdetermined_problem}, then we have a bound from above for $|\nabla \varphi|$:
\begin{equation}
\label{eq:bound_on_|nabla_varphi|_in_solutions_of_the_overdetermined_problem}
\frac{\mathcal H_{N - 1}(\Gamma_\varphi)}{N} < \int_\omega \frac{1}{\sqrt{1 + |\nabla \varphi|^2}} \ dx'.
\end{equation}
\end{proposition}

\begin{proof}
From the Minkowski formula \eqref{eqthesisMink} we have
\begin{align}
\int_{\Gamma_\varphi} H_{\Gamma_\varphi} \langle x_N e_N, \nu \rangle \ d\sigma 
&=\frac{1}{N - 1} \int_\omega \frac{|\nabla \varphi|^2}{\sqrt{1 + |\nabla \varphi|^2}} \ d x' \nonumber \\
&= \frac{1}{N - 1} \int_\omega \frac{1 + |\nabla \varphi|^2}{\sqrt{1 + |\nabla \varphi|^2}} \ dx' - \frac{1}{N - 1} \int_{\omega} \frac{1}{\sqrt{1 + |\nabla \varphi|^2}} \ d x' \nonumber \\
& = \frac{\mathcal H_{N - 1}(\Gamma_\varphi)}{N - 1} - \frac{1}{N - 1} \int_\omega \frac{1}{\sqrt{1 + |\nabla \varphi|^2}} \ dx'.
\label{eq:id_from_Mink}
\end{align}

Now, from the bound on $H_{\Gamma_\varphi}$ given by Proposition \ref{prop:bound_on_H_in_Gamma} (observe that the case \eqref{eq:const_H_from_P}, i.e. $H_{\Gamma_\varphi} \equiv \frac{1}{Nc}$, cannot happen, in view of Proposition \ref{prop:CMC_surfaces}) we have
\begin{equation} \label{eq:inequality_H_from_Mink}
\begin{array}{lll}
 \displaystyle\int_{\Gamma_\varphi} H_{\Gamma_\varphi} \langle x_N e_N, \nu\rangle \ d\sigma 
&<& \displaystyle  \int_{\Gamma_\varphi} \frac{\langle x_N e_N, \nu \rangle}{Nc} \ d \sigma \\[12pt]
& =&  \displaystyle \frac{1}{Nc} \int_\omega \varphi \ dx' \\[12pt]
& =& \displaystyle \frac{|\Omega_\varphi|}{Nc}= \frac{\mathcal H_{N - 1}(\Gamma_\varphi)}{N},
\end{array}
\end{equation}
by \eqref{eq:|Omega|_equal_c_|Gamma|}. Finally, combining \eqref{eq:id_from_Mink} and \eqref{eq:inequality_H_from_Mink}  we get
$$
\frac{\mathcal H_{N - 1}(\Gamma_\varphi)}{N - 1} - \frac{1}{N - 1} \int_\omega \frac{1}{\sqrt{1 + |\nabla \varphi|^2}} \ d x' < \frac{\mathcal H_{N - 1}(\Gamma_\varphi)}{N}, 
$$
and \eqref{eq:bound_on_|nabla_varphi|_in_solutions_of_the_overdetermined_problem} readily follows.
\end{proof}


%
%
%
%

\vskip2em

\bibliographystyle{acm}
\bibliography{Bib}

\end{document}